\documentclass[review]{elsarticle}
\usepackage{amsmath,amsfonts,amssymb,amsthm,hyperref,color}
\usepackage[mathlines]{lineno}
\newtheorem{theorem}{Theorem}[section]

\newtheorem{remark}{Remark}[section]
\newtheorem{definition}{Definition}[section]

\modulolinenumbers[1]
\journal{Journal of \LaTeX\ Templates}
\numberwithin{equation}{section}
\begin{document}

\begin{frontmatter}

\title{A free boundary problem for a class of nonlinear nonautonomous size-structured population model}


\author[mymainaddres]{Wenbin Lv}\corref{mycorrespondingauthor}
\cortext[mycorrespondingauthor]{Corresponding author}
\ead{lvwenbin@whu.edu.cn}

\author[addres2]{Shaohua Wu}
\ead{wush8@sina.com}

\address[mymainaddress]{School of Mathematical Sciences, Shanxi University, Taiyuan 030006, China}
\address[addres2]{School of Mathematics and Statistics, Wuhan University, Wuhan 430072, China}

\begin{abstract}
In this paper, we study a free boundary problem for a class of nonlinear nonautonomous size-structured population model. Using the comparison principle and upper-lower solution methods, we establish the existence of the solution for such kind of a model.
\end{abstract}

\begin{keyword}
\texttt population model\sep free boundary problem\sep upper-lower solution methods\sep existence
\MSC[2010] 35A01\sep 35D30\sep 35F20\sep 35R35
\end{keyword}
\end{frontmatter}

\section{Introduction}

Free boundary problems deal with solving partial differential equations (PDEs) in a domain, a part of whose boundary is unknown in advance; that portion of the boundary is called a free boundary. In addition to the standard boundary conditions that are needed in order to solve the PDEs, an additional condition must be imposed at the free boundary. One then seeks to determine both the free boundary and the solution of the differential equations. The theory of free boundaries has seen great progress in the last century \cite{C-S-V2015,A.F.2000,M.E.T.2013}. Recent decades have witnessed a rapid widening of the subject area by incorporation of important free boundary topics coming from different areas: In Finance, free boundary problems appear to determine the optimal exercise value in Black-Scholes models \cite{G-Y2016,Q-C-L-Y2015}; In Mathematical Biology, they indicate the moving fronts of populations or tumors \cite{C-W2012,S.B.C.2014,D-L2014}.

Formulating mathematical models that incorporate internal and external factors of population growth is one of the best way to analysising the problem of population dynamics. The literature on the population models with various level of complexity is quite vast and detailed review is beyond the scope of this paper. We mention only some of the well-established models that have been developed over the years. Among the unstructured models, the Malthus model of exponential growth and the Verhulst logistic model are especially important. For the size-structured models with density-dependency or time-dependency we refer to \cite{B-K-W1994,S-S1967}.

In this paper, we study the following free boundary problem for a class of nonlinear nonautonomous size-structured population model
\begin{equation}\label{8}
\begin{cases}
u_t+\big(V(x,t)u\big)_x=-m\Big(x,t,P\big(u(\cdot,t)\big)\Big)u,&\quad x\in[0,h(t)],~t\in[0,T],\\
\displaystyle V(0,t)u(0,t)=C(t)+\int_0^{h(t)}\beta(x,t)u(x,t)dx,&\quad t\in[0,T],\\
u(x,0)=u_0(x),&\quad x\in[0,b],\\
h'(t)=V\big(h(t),t\big),&\quad t\in[0,T],
\end{cases}
\end{equation}
where
\begin{itemize}
  \item $u(x,t)$ is an unknown function which represents the density of individuals at times $t$.
  \item $x=h(t)$ is an unknown function which represents the free boundary.
\end{itemize}
In addition, the quantities appearing in the system \eqref{8} are some given functions.
\begin{itemize}
  \item The parameters $\beta(x,t)$ and $V(x,t)$ are the time and size-dependent reproduction and growth rates, respectively.
  \item $m$ represents the mortality rate of an individual of size $x$ at time $t$ which depends on the population measure $P\big(u(\cdot,t)\big)=\int_0^{h(t)}\eta(y)u(y,t)dy$.
  \item $C(t)$ represents the inflow of $0$-size individuals from an external source (eg. seeds carried by wind).
\end{itemize}

The fixed boundary problem
\begin{equation}\label{9}
\begin{cases}
u_t+\big(V(x,t)u\big)_x=-m\Big(x,t,P\big(u(\cdot,t)\big)\Big)u,&\quad x\in[0,L],~t\in[0,T],\\
\displaystyle V(0,t)u(0,t)=C(t)+\int_0^L\beta(x,t)u(x,t)dx,&\quad t\in[0,T],\\
u(x,0)=u_0(x),&\quad x\in[0,L],
\end{cases}
\end{equation}
corresponding to \eqref{8} have been considered extensively. The above problem \eqref{9} arises in many important applications to biology and chemistry. For example, the evolution of a size-structured population, i.e., a population where individuals are distinguished by size, has been formulated into \eqref{9} \cite{S-S1967}. On the other hand, the dynamics of coagulating particles has been investigated using similar equations as in \eqref{9} \cite{A.S.A.1997}. Existence theory have been established using the characteristic method with the contraction mapping principle \cite{C-S1995}, the semigroup theorem \cite{B-K-W1994} or the upper-lower solution method \cite{A-D2000}.

Comparing to the problem, free boundary problem is more reasonable. The range of $x$ is occupied in a fixed domain in standard literatures. In other words, the biggest size of the individual is fixed.  However, the real phenomena does not obey the law. As time goes on, the biggest size of the individual maybe change due to the environment around them. For example, $x$ represents the age of the individual. If the external environment is better, the maximum age is bigger.  To the best of our knowledge, there's few results on the free boundary problem for a class of nonlinear nonautonomous size-structured population model and this is the motivation of our work. In the present paper, we will mainly clarify the problem and show the existence and uniqueness of the solution for such kind of model.

This paper is arranged as follows. In section \ref{sec2}, we present the vital conditions and the main result of the paper. In section \ref{sec3}, we take a classical transform to straighten the free boundary. A comparison result is discussed for the problem. Using upper and lower solution method, we establish the existence of the solution for such kind of the model. In addition, the uniqueness of the solution is discussed in the following section \ref{sec4}.

\section{Main result}\label{sec2}

In this section, we state our main result. Before stating our theorem, we assume some vital conditions on the parameters in \eqref{8}.

\begin{itemize}
  \item[(C1)] $V(x,t)$ is continuously differentiable with respect to $x$ and $t$. Furthermore, $V(x,t)>0$ for $(x,t)\in[0,h(t))\times[0,T]$, $V(h(t),t)=0$ and $V(\lambda h(t),t)>\lambda V(h(t),t)$ for any $0<\lambda<1$.
  \item[(C2)] $P\big(u(\cdot,t)\big)=\int_0^{h(t)}\eta(y)u(y,t)dy$ where $\eta\in L^\infty(0,h(t))$ and $\eta\geqslant0$ a.e. in $(0,h(t))$.
  \item[(C3)] $m(x,t,P)(\geqslant0)$ is continuous with respect to $x$ and $t$ and
continuously differentiable with respect to $P$ for $(x,t,P)\in[0,h(t)]\times[0,T]\times[0,+\infty)$. In addition, there exists a constant $M>0$ such that $M+m_P(x,t,P)\geqslant0$.
  \item[(C4)] $C(t)(\geqslant0)$ is continuous for $t\in[0,T]$.
  \item[(C5)] $\beta(x,t)(\geqslant0)$ is continuous with respect to $x$ and $t$ for $(x,t)\in[0,h(t)]\times[0,T]$.
\end{itemize}

\begin{remark}\label{rm1}
The main difference of the conditions on the parameters in \eqref{8} is (C2), i.e.
\[P\big(u(\cdot,t)\big)=\int_0^{h(t)}\eta(y)u(y,t)dy.\]
As we all knowm, we often assume $P\big(u(\cdot,t)\big)=\int_0^L\eta(y)u(y,t)dy$ which means the biggest size of the individual is fixed in the known results. However, the real phenomena does not obey the law. As time goes on, the biggest size of the individual maybe change due to the environment around them. For example, $x$ represents the age of the individual. If the external environment is better, the maximum age is bigger. Hence, the population measure would be (C2) and it is more reasonable compared to the known assumptions in the literature.
\end{remark}

Our main result is as follows.

\begin{theorem}\label{th2}
Assume that (C1)-(C5) hold. If
$$u_0\in L^\infty(0,b),\quad u_0\geqslant0~\text{a.e.}~\text{in}~(0,b),$$
where $0<b<L$ is a constant, then there exist an unique solution $u(x,t)$ and an unique curve $x=h(t)$ which satisfy \eqref{8}. Moreover, $x=h(t)$ is an increasing function and the solution is global in time.
\end{theorem}

\begin{remark}
The monotone of the function $x=h(t)$ is obvious, so we only to show the existence and uniqueness of the solution for the model \eqref{8}.
\end{remark}

\section{The proof of the main theorem}

In this section, we show the existence and uniqueness of the solution for a class of nonlinear nonautonomous size-structured population model \eqref{8}, i.e. Theorem \ref{th2}.

Before starting our main contents, we give simple description of the approaches. In order to achieve the goal, we shall processed as follows:
\begin{itemize}
\item First, we show the existence and uniqueness of the free boundary by the standard ODE theory;
\item Second, we straighten the free boundary and convert the problem to a fixed boundary problem. Then, we can show the existence of the solution $u(x,t)$ by comparison principle, monotone sequences and lower and upper solutions methods;
\item Finally, we obtain the solution is indeed unique.
\end{itemize}

\subsection{Existence of the solution}\label{sec3}

\subsubsection{Existence and uniqueness of the free boundary}

Noticing $h'(t)=V\big(h(t),t\big)$ and $V(x,t)$ is continuously differentiable with respect to $x$ and $t$, we get an unique continuous solution $h(t)$ by ODE standard theory \cite{H.B.2010}.

\subsubsection{Straighten the free boundary}

For that $h(t)$, we consider the problem
\begin{equation*}
\begin{cases}
u_t+\big(V(x,t)u\big)_x=-m\Big(x,t,P\big(u(\cdot,t)\big)\Big)u,&\quad x\in[0,h(t)],~t\in[0,T],\\
\displaystyle V(0,t)u(0,t)=C(t)+\int_0^{h(t)}\beta(x,t)u(x,t)dx,&\quad t\in[0,T],\\
u(x,0)=u_0(x),&\quad x\in[0,b],
\end{cases}
\end{equation*}
Take the transform
$$\xi=\frac{x}{h(t)},$$
and set
$$\widetilde{u}(\xi,t)=u\big(h(t)\xi,t\big)=u(x,t).$$
A simple calculation shows
$$\widetilde{u}_t=u_t+\xi h'(t)u_x,$$
$$\widetilde{u}_\xi=h(t)u_x.$$
Thus $\widetilde{u}$ satisfies the following equation
\begin{equation}\label{1}
\begin{cases}
\displaystyle\widetilde{u}_t+\left(\frac{V\big(h(t)\xi,t\big)-\xi h'(t)}{h(t)}\widetilde{u}\right)_\xi=-\left(m\Big(h(t)\xi,t,\widetilde{P}\big(\widetilde{u}(\cdot,t)\big)\Big)+\frac{h'(t)}{h(t)}\right)\widetilde{u},&\quad \xi\in[0,1],~t\in[0,T],\\
\displaystyle V(0,t)\widetilde{u}(0,t)=C(t)+h(t)\int_0^1\beta\big(h(t)\xi,t\big)\widetilde{u}(\xi,t)d\xi,&\quad t\in[0,T],\\
\widetilde{u}(\xi,0)=u_0(b\xi),&\quad \xi\in[0,1].
\end{cases}
\end{equation}
where $\widetilde{P}\big(\widetilde{u}(\cdot,t)\big)=h(t)\int_0^1\eta\big(h(t)\xi\big)\widetilde{u}(\xi,t)d\xi$. To solve the problem, we follow the similar argument \cite{A-D1997,A-D2000}.

\subsubsection{Comparison principle}

Let $D_T=(0,T)\times(0,1)$, and we introduce the following definition of a pair of coupled upper and lower solutions of problem \eqref{1}.

\begin{definition}
A pair of functions $\widetilde{u}(\xi,t)$ and $\widetilde{v}(\xi,t)$ are called an upper and a lower solution of \eqref{1} on $D_T$, respectively, if all the following hold:
\begin{itemize}
  \item $\widetilde{u},\widetilde{v}\in L^\infty(D_T)$;
  \item $\widetilde{u}(\xi,0)\geqslant u_0(b\xi)\geqslant\widetilde{v}(\xi,0)$ almost everywhere in $(0,1)$;
  \item For every $t\in(0,T)$ and every nonnegative $\varphi\in C^1(\overline{D_T})$, we have
  \begin{equation*}
  \begin{aligned}
   &\int_0^1\widetilde{u}(\xi,t)\varphi(\xi,t)d\xi\geqslant\int_0^1\widetilde{u}(\xi,0)\varphi(\xi,0)d\xi\\
   &+\int_0^t\left(\frac{C(s)}{h(s)}+\int_0^1\beta\big(h(s)\xi,s\big)\widetilde{u}(\xi,s)d\xi\right)\varphi(0,s)ds\\
   &+\int_0^t\int_0^1\left(\varphi_s(\xi,s)+\frac{V\big(h(s)\xi,s\big)-\xi h'(s)}{h(s)}\varphi_\xi(\xi,s)\right)\widetilde{u}(\xi,s)d\xi ds\\
   &-\int_0^t\int_0^1\left(m\Big(h(s)\xi,s,\widetilde{P}\big(\widetilde{v}(\cdot,s)\big)\Big)+M\widetilde{P}\big(\widetilde{v}(\cdot,s)\big)-M\widetilde{P}\big(\widetilde{u}(\cdot,s)\big)+\frac{h'(s)}{h(s)}\right)\widetilde{u}(\xi,s)\varphi(\xi,s)d\xi ds
  \end{aligned}
  \end{equation*}
  and
  \begin{equation*}
  \begin{aligned}
   &\int_0^1\widetilde{v}(\xi,t)\varphi(\xi,t)d\xi\leqslant\int_0^1\widetilde{v}(\xi,0)\varphi(\xi,0)d\xi\\
   &+\int_0^t\left(\frac{C(s)}{h(s)}+\int_0^1\beta\big(h(s)\xi,s\big)\widetilde{v}(\xi,s)d\xi\right)\varphi(0,s)ds\\
   &+\int_0^t\int_0^1\left(\varphi_s(\xi,s)+\frac{V\big(h(s)\xi,s\big)-\xi h'(s)}{h(s)}\varphi_\xi(\xi,s)\right)\widetilde{v}(\xi,s)d\xi ds\\
   &-\int_0^t\int_0^1\left(m\Big(h(s)\xi,s,\widetilde{P}\big(\widetilde{u}(\cdot,s)\big)\Big)+M\widetilde{P}\big(\widetilde{u}(\cdot,s)\big)-M\widetilde{P}\big(\widetilde{v}(\cdot,s)\big)+\frac{h'(s)}{h(s)}\right)\widetilde{v}(\xi,s)\varphi(\xi,s)d\xi ds.
  \end{aligned}
  \end{equation*}
\end{itemize}
\end{definition}

\begin{definition}
A function $\widetilde{u}(\xi,t)$ is called a weak solution of \eqref{1} on $D_T$ if $\widetilde{u}$ is not only an upper solution but also a lower solution of \eqref{1} on $D_T$.
\end{definition}

\begin{theorem}\label{th1}
Suppose that the assumptions in Theorem \ref{th2} hold. Let $\widetilde{u}$ and $\widetilde{v}$ be a nonnegative upper solution and a nonnegative lower solution of \eqref{1}, respectively. Then $\widetilde{u}\geqslant \widetilde{v}$ almost everywhere in $D_T$.
\end{theorem}

\begin{proof}
Let $w=\widetilde{v}-\widetilde{u}$, then $w$ satisfies
\begin{equation}\label{4}
w(\xi,0)=\widetilde{v}(\xi,0)-\widetilde{u}(\xi,0)\leqslant0\quad\text{a.e.}\quad\text{in}~(0,1),
\end{equation}
and
\begin{equation*}
  \begin{aligned}
   &\int_0^1w(\xi,t)\varphi(\xi,t)d\xi\leqslant\int_0^1w(\xi,0)\varphi(\xi,0)d\xi\\
   &+\int_0^t\left(\int_0^1\beta\big(h(s)\xi,s\big)w(\xi,s)d\xi\right)\varphi(0,s)ds\\
   &+\int_0^t\int_0^1\left(\varphi_s(\xi,s)+\frac{V\big(h(s)\xi,s\big)-\xi h'(s)}{h(s)}\varphi_\xi(\xi,s)\right)w(\xi,s)d\xi ds\\
   &-\int_0^t\int_0^1m\Big(h(s)\xi,s,\widetilde{P}\big(\widetilde{u}(\cdot,s)\big)\Big)w(\xi,s)\varphi(\xi,s)d\xi ds\\
   &+\int_0^t\int_0^1\Big(M+m_{\widetilde{P}}\big(\xi,s,\theta(s)\big)\Big)\widetilde{P}\big(w(\cdot,s)\big)\widetilde{u}(\xi,s)\varphi(\xi,s)d\xi ds\\
   &+\int_0^t\int_0^1M\widetilde{P}\big(w(\cdot,s)\big)\widetilde{v}(\xi,s)\varphi(\xi,s)d\xi ds-\int_0^t\int_0^1\frac{h'(s)}{h(s)}w(\xi,s)\varphi(\xi,s)d\xi ds,
  \end{aligned}
\end{equation*}
by mean value theorem, where $\theta(s)$ between $\widetilde{P}\big(\widetilde{u}(\cdot,s)\big)$ and $\widetilde{P}\big(\widetilde{v}(\cdot,s)\big)$.

Let $\varphi(\xi,t)=e^{\lambda t}\psi(\xi,t)$ where $\psi\in C^1(\overline{D_T})$ and
$\lambda(>0)$ is chosen so that $\lambda-m-\frac{h'(t)}{h(t)}\geqslant0$ on $D_T\times[\widetilde{P}_1,\widetilde{P}_2]$ where $$\widetilde{P}_1=\min\{\inf_{[0,T]}\widetilde{P}(\widetilde{u}(\cdot,t)),\inf_{[0,T]}\widetilde{P}(\widetilde{v}(\cdot,t))\}$$ and $$\widetilde{P}_2=\max\{\sup_{[0,T]}\widetilde{P}(\widetilde{u}(\cdot,t)),\sup_{[0,T]}\widetilde{P}(\widetilde{v}(\cdot,t))\}.$$ Then we find
\begin{equation}\label{3}
 \begin{aligned}
   &e^{\lambda t}\int_0^1w(\xi,t)\psi(\xi,t)d\xi\leqslant\int_0^1w(\xi,0)\psi(\xi,0)d\xi\\
   &+\int_0^t\left(\int_0^1\beta\big(h(s)\xi,s\big)w(\xi,s)d\xi\right)e^{\lambda s}\psi(0,s)ds\\
   &+\int_0^t\int_0^1\left(\psi_s(\xi,s)+\frac{V\big(h(s)\xi,s\big)-\xi h'(s)}{h(s)}\psi_\xi(\xi,s)\right)e^{\lambda s}w(\xi,s)d\xi ds\\
   &+\int_0^t\int_0^1\left(\lambda-m\Big(h(s)\xi,s,\widetilde{P}\big(\widetilde{u}(\cdot,s)\big)\Big)-\frac{h'(s)}{h(s)}\right)w(\xi,s)e^{\lambda s}\psi(\xi,s)d\xi ds\\
   &+\int_0^t\int_0^1\Big(M+m_{\widetilde{P}}\big(\xi,s,\theta(s)\big)\Big)\widetilde{P}\big(w(\cdot,s)\big)\widetilde{u}(\xi,s)e^{\lambda s}\psi(\xi,s)d\xi ds\\
   &+\int_0^t\int_0^1M\widetilde{P}\big(w(\cdot,s)\big)\widetilde{v}(\xi,s)e^{\lambda s}\psi(\xi,s)d\xi ds.
 \end{aligned}
\end{equation}
To simplify the above inequality, we now set up a backward problem as follows:
\begin{equation}\label{2}
\begin{cases}
\psi_s(\xi,s)+\frac{V\big(h(s)\xi,s\big)-\xi h'(s)}{h(s)}\psi_\xi(\xi,s)=0,&0<s<t,~0<\xi<1,\\
\psi(1,s)=0,&0<s<t,\\
\psi(\xi,t)=\zeta(\xi),&0\leqslant\xi\leqslant1.
\end{cases}
\end{equation}
Here $\zeta(\xi)\in C^\infty_0(0,1),0\leqslant\zeta\leqslant1$. In order to solve the problem \eqref{2}, we take a transform $\tau=t-s$ and let $\widetilde{\psi}(\xi,\tau)=\psi(\xi,s)$. Then \eqref{2} can be written as
\begin{equation}\label{5}
\begin{cases}
\widetilde{\psi}_\tau(\xi,\tau)-\frac{V\big(h(t-\tau)\xi,t-\tau\big)+\xi h'(t-\tau)}{h(t-\tau)}\widetilde{\psi}_\xi(\xi,\tau)=0,&0<\tau<t,~0<\xi<1,\\
\widetilde{\psi}(1,\tau)=0,&0<\tau<t,\\
\widetilde{\psi}(\xi,0)=\zeta(\xi),&0\leqslant\xi\leqslant1.
\end{cases}
\end{equation}
\eqref{5} can be solved by the characteristic method. Note that the initial and boundary values for $\zeta$ imply that $0\leqslant\psi\leqslant1$ on $D_T$.

Substituting such a $\psi$ in \eqref{3} yields
$$\int_0^1w(\xi,t)\zeta(\xi)d\xi\leqslant\int_0^1w(\xi,0)^+d\xi+C_1\int_0^t\int_0^1w(\xi,s)^+d\xi ds,$$
where $w(\xi,t)^+=\max\{w(\xi,t),0\}$ and $$C_1=\max_{D_T}\left\{\beta\big(h(t)\xi,t\big)+\left(\lambda-m-\frac{h'(t)}{h(t)}\right)+\left[(M+m_{\widetilde{P}})\widetilde{u}(\xi,t)+M\widetilde{v}(\xi,t)\right]h(t)\eta\big(h(t)\xi\big)\right\}.$$
From the condition on
initial data in \eqref{4}, we have
$$\int_0^1w(\xi,t)\zeta(\xi)d\xi\leqslant C_1\int_0^t\int_0^1w(\xi,s)^+d\xi ds.$$
Since this inequality holds for every $\zeta\in C_0^\infty(0,1)$ with $0\leqslant\zeta\leqslant1$,
we can choose a sequence $\{\zeta_n\}_{n=1}^\infty$ on $(0,1)$ converging to
$$\chi=\begin{cases}
1,&w(\xi,t)>0,\\
0,&\text{otherwise}.
\end{cases}$$
Consequently, we find that
$$\int_0^1w(\xi,t)^+d\xi\leqslant C_1\int_0^t\int_0^1w(\xi,s)^+d\xi ds,$$
which by Gronwall¡¯s inequality leads to
$$\int_0^1w(\xi,t)^+d\xi=0.$$
Thus, the proof is completed.
\end{proof}

\subsubsection{Monotone sequences and existence of solutions}

Now, we construct a pair of nonnegative lower and upper solutions of \eqref{1}.

Let $\underline{u}^0(\xi,t)=0$ and $\overline{u}^0(\xi,t)=\delta e^{\sigma t}e^{-\gamma \xi}$, where $\delta,\sigma,\gamma$ are some determinate constants. Then it can be easily shown that $\underline{u}^0$ and $\overline{u}^0$ are a pair of coupled lower and upper solutions of \eqref{1} on $[0,1]\times[0,T_0]$ with $T_0=\min\{T,\frac{\ln 2}{\sigma}\}$.

Actually, it is easily seen that $\underline{u}^0$ is a lower solution of \eqref{1}. The task is now to show that $\overline{u}^0$ is an upper solution of \eqref{1}.
\begin{enumerate}
  \item $\overline{u}^0\in L^\infty(D_T)$.
  \item $u_0(b\xi)\leqslant\delta e^{-\gamma}\leqslant\overline{u}^0(\xi,0)$ according to the choice of the parameters $\delta$ and $\gamma$ in the following.
  \item Notice
  $$\overline{u}^0_t(\xi,t)=\sigma\overline{u}^0(\xi,t),\quad\overline{u}^0_\xi(\xi,t)=-\gamma\overline{u}^0(\xi,t).$$
  We have
  \begin{align*}
    &\int_0^t\left(\frac{C(s)}{h(s)}+\int_0^1\beta\big(h(s)\xi,s\big)\overline{u}^0(\xi,s)d\xi\right)\varphi(0,s)ds\\
    =&\int_0^t\frac{C(s)}{\delta e^{\sigma s}V(0,s)}\frac{V(0,s)}{h(s)}\overline{u}^0(0,s)\varphi(0,s)ds\\
    &\quad+\int_0^t\int_0^1\frac{\beta\big(h(s)\xi,s\big)h(s)}{V(0,s)}e^{-\gamma\xi}d\xi\frac{V(0,s)}{h(s)}\overline{u}^0(0,s)\varphi(0,s)ds.
  \end{align*}
  Choose a constant $\gamma>0$ large enough such that
  $$\max_{\overline{D_T}}\frac{\beta\big(h(s)\xi,s\big)h(s)}{V(0,s)}\leqslant\frac{\gamma}{2},$$
  and then choose $\delta>0$ large enough such that
  $$\|u_0\|_{L^\infty}\leqslant\delta e^{-\gamma} \quad \text{and} \quad \max_{[0,T]}\frac{C(s)}{V(0,s)}\leqslant\frac{\delta}{2}.$$
  Hence, it holds that
  $$\int_0^t\left(\frac{C(s)}{h(s)}+\int_0^1\beta\big(h(s)\xi,s\big)\overline{u}^0(\xi,s)d\xi\right)\varphi(0,s)ds\leqslant\int_0^t\frac{V(0,s)}{h(s)}\overline{u}^0(0,s)\varphi(0,s)ds.$$
  We have
  \begin{align*}
    &\int_0^t\int_0^1\varphi_s(\xi,s)\overline{u}^0(\xi,s)d\xi ds\\
    =&\int_0^1\varphi(\xi,t)\overline{u}^0(\xi,t)d\xi-\int_0^1\varphi(\xi,0)\overline{u}^0(\xi,0)d\xi\\
    &\quad\quad-\sigma\int_0^t\int_0^1\varphi(\xi,s)\overline{u}^0(\xi,s)d\xi ds,
  \end{align*}
  and
  \begin{align*}
    &\int_0^t\int_0^1\frac{V\big(h(s)\xi,s\big)-\xi h'(s)}{h(s)}\varphi_\xi(\xi,s)\overline{u}^0(\xi,s)d\xi ds\\
    =&\int_0^t\frac{V\big(h(s),s\big)-h'(s)}{h(s)}\varphi(1,s)\overline{u}^0(1,s)ds-\int_0^t\frac{V(0,s)}{h(s)}\varphi(0,s)\overline{u}^0(0,s)ds\\
    &\quad\quad-\int_0^t\int_0^1\frac{V_x\big(h(s)\xi,s\big)h(s)-h'(s)}{h(s)}\varphi(\xi,s)\overline{u}^0(\xi,s)d\xi ds\\
    &\quad\quad+\gamma\int_0^t\int_0^1\frac{V\big(h(s)\xi,s\big)-\xi h'(s)}{h(s)}\varphi(\xi,s)\overline{u}^0(\xi,s)d\xi ds.
  \end{align*}
  by integrating by parts. Then, it holds that
  \begin{equation*}
  \begin{aligned}
   &\int_0^1\overline{u}^0(\xi,0)\varphi(\xi,0)d\xi+\int_0^t\left(\frac{C(s)}{h(s)}+\int_0^1\beta\big(h(s)\xi,s\big)\overline{u}^0(\xi,s)d\xi\right)\varphi(0,s)ds\\
   &+\int_0^t\int_0^1\left(\varphi_s(\xi,s)+\frac{V\big(h(s)\xi,s\big)-\xi h'(s)}{h(s)}\varphi_\xi(\xi,s)\right)\overline{u}^0(\xi,s)d\xi ds\\
   &-\int_0^t\int_0^1\left(m\Big(h(s)\xi,s,\widetilde{P}\big(\underline{u}^0(\cdot,s)\big)\Big)+M\widetilde{P}\big(\underline{u}^0(\cdot,s)\big)-M\widetilde{P}\big(\overline{u}^0(\cdot,s)\big)+\frac{h'(s)}{h(s)}\right)\overline{u}^0(\xi,s)\varphi(\xi,s)d\xi ds\\
   \leqslant&\int_0^1\overline{u}^0(\xi,t)\varphi(\xi,t)d\xi\\
   &+\int_0^t\int_0^1\left(-\sigma-V_x\big(h(s)\xi,s\big)+\gamma\frac{V\big(h(s)\xi,s\big)-\xi h'(s)}{h(s)}+Mh(s)\int_0^1\eta(y)\overline{u}^0(y,s)dy\right)\varphi(\xi,s)\overline{u}^0(\xi,s)d\xi ds,
  \end{aligned}
  \end{equation*}
  by the fourth equation of \eqref{8}. Choose $\sigma>0$ large enough such that
  $$\sigma\geqslant\max_{\overline{D_T}}\left|V_x\big(h(s)\xi,s\big)\right|+\gamma\max_{\overline{D_T}}\frac{V\big(h(s)\xi,s\big)-\xi h'(s)}{h(s)}+2\delta M\|\eta\|_{L^\infty}\max_{[0,T]}h(s).$$
\end{enumerate}
Then $\overline{u}^0$ is an upper solution of \eqref{1}.

We then define two sequences $\{\underline{u}^k\}^\infty_{k=0}$ and $\{\overline{u}^k\}^\infty_{k=0}$ as follows.

For $k=1,2,\cdots$,
\begin{equation}\label{6}
\begin{cases}
\displaystyle\underline{u}^k_t+\left(\frac{V\big(h(t)\xi,t\big)-\xi h'(t)}{h(t)}\underline{u}^k\right)_\xi=-\left(A^{k-1}(\xi,t)+\frac{h'(t)}{h(t)}\right)\underline{u}^k,&\quad \xi\in[0,1],~t\in[0,T],\\
\displaystyle V(0,t)\underline{u}^k(0,t)=C(t)+h(t)\int_0^1\beta\big(h(t)\xi,t\big)\underline{u}^{k-1}(\xi,t)d\xi,&\quad t\in[0,T],\\
\underline{u}^k(\xi,0)=u_0(b\xi),&\quad \xi\in[0,1],
\end{cases}
\end{equation}
and
\begin{equation}\label{7}
\begin{cases}
\displaystyle\overline{u}^k_t+\left(\frac{V\big(h(t)\xi,t\big)-\xi h'(t)}{h(t)}\overline{u}^k\right)_\xi=-\left(B^{k-1}(\xi,t)+\frac{h'(t)}{h(t)}\right)\overline{u}^k,&\quad \xi\in[0,1],~t\in[0,T],\\
\displaystyle V(0,t)\overline{u}^k(0,t)=C(t)+h(t)\int_0^1\beta\big(h(t)\xi,t\big)\overline{u}^{k-1}(\xi,t)d\xi,&\quad t\in[0,T],\\
\overline{u}^k(\xi,0)=u_0(b\xi),&\quad \xi\in[0,1],
\end{cases}
\end{equation}
where
$$A^{k-1}(\xi,t)=\left(m\Big(h(t)\xi,t,\widetilde{P}\big(\overline{u}^{k-1}(\cdot,t)\big)\Big)\right)-M\widetilde{P}\big(\underline{u}^{k-1}(\cdot,s)\big)+M\widetilde{P}\big(\overline{u}^{k-1}(\cdot,s)\big),$$
$$B^{k-1}(\xi,t)=\left(m\Big(h(t)\xi,t,\widetilde{P}\big(\underline{u}^{k-1}(\cdot,t)\big)\Big)\right)-M\widetilde{P}\big(\overline{u}^{k-1}(\cdot,s)\big)+M\widetilde{P}\big(\underline{u}^{k-1}(\cdot,s)\big).$$
The existence of solutions for the problems \eqref{6} and \eqref{7} follows from the method of characteristics.  Consider the equation for the characteristic curves given by
\begin{equation*}
\begin{cases}
\frac{d}{ds}t(s)=1,\\
\frac{d}{ds}\xi(s)=\frac{V\Big(h\big(t(s)\big)\xi(s),t(s)\Big)-\xi(s)h'\big(t(s)\big)}{h\big(t(s)\big)}.
\end{cases}
\end{equation*}
The solution $\underline{u}^k$ of \eqref{6} along the characteristic curve $\big(\xi(s),t(s)\big)$ satisfies the following equation
$$\frac{d}{ds}\underline{u}^k(s)=-\left(V_x\Big(h\big(t(s)\big)\xi(s),t(s)\Big)+A^{k-1}\big(\xi(s),t(s)\big)\right)\underline{u}^k(s).$$
Parametrizing the characteristic curves with the variable $t$, then a characteristic curve passing through $(\widehat{\xi},\widehat{t})$ is given by $(t,X(t;\widehat{\xi},\widehat{t}))$ where $X$ satisfies
$$\frac{d}{dt}X(t;\widehat{\xi},\widehat{t})=\frac{V\big(h(t)X(t;\widehat{\xi},\widehat{t}),t\big)-X(t;\widehat{\xi},\widehat{t})h'(t)}{h(t)}$$
and $X(\widehat{t};\widehat{\xi},\widehat{t})=\widehat{\xi}$. From (C1), it follows that the function $X$ is strictly increasing. Hence, a unique inverse function $\tau(\xi;\widehat{t},\widehat{\xi})$ exists. Now we define $G(\xi)=\tau(\xi;0,0)$ where $\big(G(\xi),\xi\big)$ represents the characteristic curve passing through $(0,0)$ which divides the $(\xi,t)$-plane into two parts. Then for any point $(\xi,t)$ with $t\leqslant G(\xi)$, the solution $\underline{u}^k(\xi,t)$ is determined through the initial condition by
\begin{equation*}
\underline{u}^k(\xi,t)=u_0\big(bX(0;\xi,t)\big)e^{-\int_0^t\left(V_x\big(h(s)X(s;\xi,t),s\big)+A^{k-1}\big(X(s;\xi,t),s\big)\right)ds}
\end{equation*}
and for any point $(\xi,t)$ with $t>G(\xi)$ the solution is determined via the boundary condition by
\begin{equation*}
\underline{u}^k(\xi,t)=R^{k-1}\big(\tau(0;\xi,t)\big)e^{-\int_{\tau(0;\xi,t)}^t\left(V_x\big(h(s)X(s;\xi,t),s\big)+A^{k-1}\big(X(s;\xi,t),s\big)\right)ds},
\end{equation*}
where $R^{k-1}(t)=\frac{C(t)+h(t)\int_0^1\beta\big(h(t)\xi,t\big)\underline{u}^{k-1}(\xi,t)d\xi}{V(0,t)}$. Similarly, for any point $(\xi,t)$ with $t\leqslant G(\xi)$, the solution $\overline{u}^k(\xi,t)$ is determined through the initial condition by
\begin{equation*}
\overline{u}^k(\xi,t)=u_0\big(bX(0;\xi,t)\big)e^{-\int_0^t\left(V_x\big(h(s)X(s;\xi,t),s\big)+B^{k-1}\big(X(s;\xi,t),s\big)\right)ds}
\end{equation*}
and for any point $(\xi,t)$ with $t>G(\xi)$ the solution is determined via the boundary condition by
\begin{equation*}
\overline{u}^k(\xi,t)=Q^{k-1}\big(\tau(0;\xi,t)\big)e^{-\int_{\tau(0;\xi,t)}^t\left(V_x\big(h(s)X(s;\xi,t),s\big)+B^{k-1}\big(X(s;\xi,t),s\big)\right)ds},
\end{equation*}
where $Q^{k-1}(t)=\frac{C(t)+h(t)\int_0^1\beta\big(h(t)\xi,t\big)\overline{u}^{k-1}(\xi,t)d\xi}{V(0,t)}$.

Next we show that the sequences $\{\underline{u}^k\}^\infty_{k=0}$ and $\{\overline{u}^k\}^\infty_{k=0}$ are monotone by induction. The procedure of induction is as follows.

Step 1: Initial hypothesis of the induction;

We first let $w=\underline{u}^0-\underline{u}^1$. Then $w$ satisfies
\begin{equation*}
  \begin{aligned}
   &\int_0^1w(\xi,t)\varphi(\xi,t)d\xi\leqslant\int_0^1w(\xi,0)\varphi(\xi,0)d\xi\\
   &+\int_0^t\left(\int_0^1\beta\big(h(s)\xi,s\big)w(\xi,s)d\xi\right)\varphi(0,s)ds\\
   &+\int_0^t\int_0^1\left(\varphi_s(\xi,s)+\frac{V\big(h(s)\xi,s\big)-\xi h'(s)}{h(s)}\varphi_\xi(\xi,s)\right)w(\xi,s)d\xi ds\\
   &-\int_0^t\int_0^1A^0(\xi,s)w(\xi,s)\varphi(\xi,s)d\xi ds-\int_0^t\int_0^1\frac{h'(s)}{h(s)}w(\xi,s)\varphi(\xi,s)d\xi ds.
  \end{aligned}
\end{equation*}
Similar to the proof of Theorem \ref{th1}, we can get $w\leqslant0$ which implies $\underline{u}^0\leqslant\underline{u}^1$. Similarly, it can be seen that $\overline{u}^0\geqslant\overline{u}^1$. From this and the facts that $\underline{u}^1$ and $\overline{u}^1$ are a lower and an upper solution, respectively, we obtain $\underline{u}^1\leqslant\overline{u}^1$.

Step 2: Hypothesis and claim of the induction;

Assume that for some $k>1$, $\underline{u}^k$ and $\overline{u}^k$ are a lower and an upper solution of \eqref{1}, respectively. By similar reasoning, we can show that $\underline{u}^k\leqslant\underline{u}^{k+1}\leqslant\overline{u}^{k+1}\leqslant\overline{u}^k$ and that $\underline{u}^{k+1}$ and $\overline{u}^{k+1}$ are also a lower and an upper solution of \eqref{1}, respectively.

Thus, by induction, we obtain two monotone sequences that satisfy
$$\underline{u}^0\leqslant\underline{u}^1\leqslant\cdots\leqslant\underline{u}^k\leqslant\overline{u}^k\leqslant\cdots\overline{u}^1\leqslant\overline{u}^0~\text{a.e.}~\text{in}~\overline{D_{T_0}}.$$
for each $k=0,1,2,\cdots$. Hence it follows from the monotonicity of the sequences $\{\underline{u}^k\}^\infty_{k=0}$ and $\{\overline{u}^k\}^\infty_{k=0}$, there exist functions $\underline{u}$ and $\overline{u}$ such that $\underline{u}^k\rightarrow\underline{u}$ and $\overline{u}^k\rightarrow\overline{u}$ pointwise in $D_{T_0}$. Clearly $\underline{u}\leqslant\overline{u}$ almost everywhere in $D_{T_0}$.

Upon establishing the monotonicity of the sequences, we now prove the sequences
$\{\underline{u}^k\}^\infty_{k=0}$ and $\{\overline{u}^k\}^\infty_{k=0}$ converge uniformly along characteristic curves to a limit function $\widetilde{u}$.

On one hand, owing to $\underline{u}^0\leqslant\underline{u}^k\leqslant\overline{u}^0$ and $\underline{u}^k\leqslant\underline{u}^{k+1}$, we obtain that, along the characteristic curve passing through $(\xi_0,0)$, the solution
$$\underline{u}^k(X(t;\xi_0,0),t)=u_0\big(b\xi_0\big)e^{-\int_0^t\left(V_x\big(h(s)X(s;\xi_0,0),s\big)+A^{k-1}\big(X(s;\xi_0,0),s\big)\right)ds}$$
converges to
$$\underline{u}(X(t;\xi_0,0),t)=u_0\big(b\xi_0\big)e^{-\int_0^t\left(V_x\big(h(s)X(s;\xi_0,0),s\big)+A\big(X(s;\xi_0,0),s\big)\right)ds}$$
uniformly and monotonically for $0\leqslant t\leqslant T_0$, where
$$A\big(X(s;\xi_0,0),s\big)=\left(m\Big(h(t)\xi,t,\widetilde{P}\big(\overline{u}(\cdot,t)\big)\Big)\right)-M\widetilde{P}\big(\underline{u}(\cdot,s)\big)+M\widetilde{P}\big(\overline{u}(\cdot,s)\big).$$
On the other hand, along the characteristic curve passing through $(0,t_0)$, the solution
$$\underline{u}^k(X(t;0,t_0),t)=R^{k-1}\big(t_0\big)e^{-\int_{t_0}^t\left(V_x\big(h(s)X(s;0,t_0),s\big)+A^{k-1}\big(X(s;0,t_0),s\big)\right)ds}$$
converges to
$$\underline{u}(X(t;0,t_0),t)=R\big(t_0\big)e^{-\int_{t_0}^t\left(V_x\big(h(s)X(s;0,t_0),s\big)+A\big(X(s;0,t_0),s\big)\right)ds}$$
uniformly and monotonically for $0\leqslant t\leqslant T_0$, where
$$R(t)=\frac{C(t)+h(t)\int_0^1\beta\big(h(t)\xi,t\big)\underline{u}(\xi,t)d\xi}{V(0,t)}.$$
Thus, we have
\begin{equation*}
\underline{u}(\xi,t)=\left\{
                       \begin{array}{ll}
                        u_0\big(bX(0;\xi,t)\big)e^{-\int_0^t\left(V_x\big(h(s)X(s;\xi,t),s\big)+A\big(X(s;\xi,t),s\big)\right)ds}, & \hbox{$t\leqslant G(\xi)$;} \\
                        R\big(\tau(0;\xi,t)\big)e^{-\int_{\tau(0;\xi,t)}^t\left(V_x\big(h(s)X(s;\xi,t),s\big)+A\big(X(s;\xi,t),s\big)\right)ds}, & \hbox{$t>G(\xi)$.}
                       \end{array}
                     \right.
\end{equation*}
Similarly, we have
\begin{equation*}
\overline{u}(\xi,t)=\left\{
                       \begin{array}{ll}
                        u_0\big(bX(0;\xi,t)\big)e^{-\int_0^t\left(V_x\big(h(s)X(s;\xi,t),s\big)+B\big(X(s;\xi,t),s\big)\right)ds}, & \hbox{$t\leqslant G(\xi)$;} \\
                        Q\big(\tau(0;\xi,t)\big)e^{-\int_{\tau(0;\xi,t)}^t\left(V_x\big(h(s)X(s;\xi,t),s\big)+B\big(X(s;\xi,t),s\big)\right)ds}, & \hbox{$t>G(\xi)$.}
                       \end{array}
                     \right.
\end{equation*}
where
\begin{gather*}
B(\xi,t)=\left(m\Big(h(t)\xi,t,\widetilde{P}\big(\underline{u}(\cdot,t)\big)\Big)\right)-M\widetilde{P}\big(\overline{u}(\cdot,s)\big)+M\widetilde{P}\big(\underline{u}(\cdot,s)\big),\\
Q(t)=\frac{C(t)+h(t)\int_0^1\beta\big(h(t)\xi,t\big)\overline{u}(\xi,t)d\xi}{V(0,t)}.
\end{gather*}

We now show that $\underline{u}(\xi,t)=\overline{u}(\xi,t)$. Let $w=\overline{u}(\xi,t)-\underline{u}(\xi,t)$. Since $\overline{u}(\xi,t)\geqslant\underline{u}(\xi,t)$, $w(\xi,t)\geqslant0$ and $w(\xi,0)=0$. Hence $w$ satisfies
\begin{equation*}
  \begin{aligned}
   \int_0^1w(\xi,t)d\xi=&\int_0^t\left(\int_0^1\beta\big(h(s)\xi,s\big)w(\xi,s)d\xi\right)ds\\
   &-\int_0^t\int_0^1m\Big(h(s)\xi,s,\widetilde{P}\big(\overline{u}(\cdot,s)\big)\Big)w(\xi,s)d\xi ds\\
   &+\int_0^t\int_0^1\Big(M+m_{\widetilde{P}}\big(\xi,s,\theta(s)\big)\Big)\widetilde{P}\big(w(\cdot,s)\big)\overline{u}(\xi,s)d\xi ds\\
   &+\int_0^t\int_0^1M\widetilde{P}\big(w(\cdot,s)\big)\underline{u}(\xi,s)d\xi ds-\int_0^t\int_0^1\frac{h'(s)}{h(s)}w(\xi,s)d\xi ds,\\
\leqslant& C_0\int_0^t\int_0^1w(\xi,s)d\xi ds
  \end{aligned}
\end{equation*}
by choosing $\varphi(\xi,t)=1$, where
$$C_0=\max_{\overline{D_{T_0}}}\left[\beta\big(h(t)\xi,t\big)+\left(M+m_{\widetilde{P}}\big(\xi,t,\theta(t)\big)\right)h(t)\|\eta\|_{L^\infty}\overline{u}(\xi,t)+Mh(t)\|\eta\|_{L^\infty}\underline{u}(\xi,t)\right].$$
Owing to Gronwall¡¯s inequality, we conclude $w(\xi,t)=0$ a.e. in $D_{T_0}$.

\subsection{Uniqueness of the solution}\label{sec4}

Using arguments in \cite{C-S1995}, we can establish $P(\widetilde{u}(\cdot,t))$ is continuous. We claim $\widetilde{u}$ is unique.

Actually, we suppose $\widetilde{u}_1(\xi,t)$ and $\widetilde{u}_2(\xi,t)$ are two nonnegative solution of \eqref{1}. If $P(\widetilde{u}_1(\cdot,t))=P(\widetilde{u}_2(\cdot,t))$ for $0<t<T$, then $\widetilde{u}_1(\xi,t)=\widetilde{u}_2(\xi,t)$. Without loss of generality, we assume that
$$\left\{
  \begin{array}{ll}
    P(\widetilde{u}_1(\cdot,t))=P(\widetilde{u}_2(\cdot,t)), & \hbox{for~$0\leqslant t\leqslant t_0$,} \\
    P(\widetilde{u}_1(\cdot,t))>P(\widetilde{u}_2(\cdot,t)), & \hbox{for~$t_0<t\leqslant t_1$,}
  \end{array}
\right.$$
where $0\leqslant t_0<t_1\leqslant T$. We have
\begin{align*}
-m(h(t)\xi,t,P(\widetilde{u}_1(\cdot,t)))=&-[m(h(t)\xi,t,P(\widetilde{u}_1(\cdot,t)))+MP(\widetilde{u}_1(\cdot,t))]+MP(\widetilde{u}_1(\cdot,t))\\
\leqslant&-m(h(t)\xi,t,P(\widetilde{u}_2(\cdot,t)))-MP(\widetilde{u}_2(\cdot,t))+MP(\widetilde{u}_1(\cdot,t))
\end{align*}
and
\begin{align*}
-m(h(t)\xi,t,P(\widetilde{u}_2(\cdot,t)))\geqslant-m(h(t)\xi,t,P(\widetilde{u}_1(\cdot,t)))-MP(\widetilde{u}_1(\cdot,t))+MP(\widetilde{u}_2(\cdot,t))
\end{align*}
by (C3). Then $\widetilde{u}_1(\xi,t)$ and $\widetilde{u}_2(\xi,t)$ are a lower and an upper solution of \eqref{1} on $D_{t_1}$, respectively. By comparison principle, we get $\widetilde{u}_1(\xi,t)<\widetilde{u}_2(\xi,t)$ a.e. in $D_{t_1}$ and hence $P(\widetilde{u}_1(\cdot,t))\leqslant P(\widetilde{u}_2(\cdot,t))$ for $0\leqslant t\leqslant t_1$, which is a contradiction.

\section*{Acknowledgments}
\addcontentsline{toc}{chapter}{Acknowledgments}
The authors of this paper would like to thank the referee for the comments and helpful suggestions.

%\bibliography{D:/paper-paper/bibfile/mybibfile}
\end{document}